\documentclass[11pt]{amsart}
\usepackage{latexsym}
\usepackage{amsmath}
\usepackage{amssymb}

\newcommand{\eproof}{\mbox{\ }\hfill $\Box$ \par \vskip 10pt}

\newtheorem{Theorem}{Theorem}[section]
\newtheorem{lemma}[Theorem]{Lemma}

\numberwithin{equation}{section}

\def\cal{\mathcal}

\begin{document}

\title[Semi-classical resolvent estimates]{Semi-classical resolvent estimates for short-range $L^\infty$ potentials}

\author[G. Vodev]{Georgi Vodev}

\address {Universit\'e de Nantes, Laboratoire de Math\'ematiques Jean Leray, 2 rue de la Houssini\`ere, BP 92208, 44322 Nantes Cedex 03, France}
\email{Georgi.Vodev@univ-nantes.fr}

\date{}

\begin{abstract} We prove semi-classical resolvent estimates for real-valued potentials $V\in L^\infty(\mathbb{R}^n)$, $n\ge 3$, 
satisfying $V(x)={\cal O}\left(\langle x\rangle^{-\delta}\right)$ with $\delta>3$. 
\end{abstract} 

\maketitle

\setcounter{section}{0}
\section{Introduction and statement of results}

Our goal in this note is to study the resolvent of the Schr\"odinger operator
$$P(h)=-h^2\Delta+V(x)$$
where $0<h\ll 1$ is a semi-classical parameter, $\Delta$ is the negative Laplacian in 
$\mathbb{R}^n$, $n\ge 3$, and $V\in L^\infty(\mathbb{R}^n)$ is a real-valued potential satisfying
\begin{equation}\label{eq:1.1}
|V(x)|\le C\langle x\rangle^{-\delta}
\end{equation}
with some constants $C>0$ and $\delta>3$. More precisely, we are interested in bounding from above the quantity
$$g_s^\pm(h,\varepsilon):=\log\left\|\langle x\rangle^{-s}(P(h)-E\pm i\varepsilon)^{-1}\langle x\rangle^{-s}
\right\|_{L^2\to L^2}$$
where $L^2:=L^2(\mathbb{R}^n)$, $0<\varepsilon<1$, $s>1/2$ and $E>0$ is a fixed energy level independent of $h$.
Such bounds are known in verious situations. For example, for long-range real-valued $C^1$ potentials it is proved in
\cite{kn:D} when $n\ge 3$ and in \cite{kn:S1} when $n=2$ that
\begin{equation}\label{eq:1.2}
g_s^\pm(h,\varepsilon)\le Ch^{-1}
\end{equation}
with some constant $C>0$ independent of $h$ and $\varepsilon$. Previously, the bound (\ref{eq:1.2}) was proved for smooth potentials
in \cite{kn:B2} and an analog of (\ref{eq:1.2}) for H\"older potentials was proved in \cite{kn:V1}. A high-frequency analog of (\ref{eq:1.2}) on more complex Riemannian manifolds was also proved in 
\cite{kn:B1} and \cite{kn:CV}. In all these papers the regularity of the potential (and of the perturbation in general) plays an essential role. 
Without any regularity the problem of bounding $g_s^\pm$ from above by an explicit function of $h$ gets quite tough. Nevertheless, it has been
recently shown in \cite{kn:S2} that for real-valued compactly supported $L^\infty$ potentials one has the bound
\begin{equation}\label{eq:1.3}
g_s^\pm(h,\varepsilon)\le Ch^{-4/3}\log(h^{-1})
\end{equation}
with some constant $C>0$ independent of $h$ and $\varepsilon$. The bound (\ref{eq:1.3}) has been also proved in \cite{kn:KV} still for
real-valued compactly supported $L^\infty$ potentials but with the weight $\langle x\rangle^{-s}$ replaced by a cut-off function. 
When $n=1$ it was shown in \cite{kn:DZ} that we have the better bound (\ref{eq:1.2}) instead of (\ref{eq:1.3}). 
When $n\ge 2$, however, the bound (\ref{eq:1.3}) seems hard to improve without extra conditions on the potential.
The problem of showing that the bound (\ref{eq:1.3}) is optimal is largely open. In contrast, it is well-known
that the bound (\ref{eq:1.2}) cannot be improved in general (e.g. see \cite{kn:DDZ}). 

In this note we show that the bound (\ref{eq:1.3}) still holds for non-compactly supported $L^\infty$ potentials
when $n\ge 3$. Our main result is the following

\begin{Theorem} Under the condition (\ref{eq:1.1}), there exists $h_0>0$ such that for all $0<h\le h_0$ 
the bound (\ref{eq:1.3}) holds true.
\end{Theorem}

\noindent
{\bf Remark.} It is easy to see from the proof (see the inequality (\ref{eq:4.2})) that the bound (\ref{eq:1.3}) holds also for a 
complex-valued potential $V$ satisfying (\ref{eq:1.1}), provided
that its imaginary part satisfies the condition
$$\mp{\rm Im}\,V(x)\ge 0\quad\mbox{for all}\quad x\in \mathbb{R}^n.$$

To prove this theorem we adapt the Carleman estimates proved in \cite{kn:S2} simplifying some key arguments 
as for example the construction of the phase function $\varphi$. This is made possible by defining the key
function $F$ in Section 3 differently, without involving the second derivative $\varphi''$. The consequence is
that we do not need to seek $\varphi'$ as a solution to a differential equation as done in \cite{kn:S2}, but 
it suffices to
define it explicitly. Note also that similar (but simpler) Carleman estimates have been used in \cite{kn:V2}
to prove high-frequency resolvent estimates for the magnetic Schr\"odinger operator with large
$L^\infty$ magnetic potentials.

\section{Construction of the phase and weight functions} 

We will first construct the weight function. We begin by introducing the continuous function
$$\mu(r)=
\left\{
\begin{array}{lll}
 (r+1)^2-1&\mbox{for}& 0\le r\le a,\\
 (a+1)^2-1+(a+1)^{-2s+1}-(r+1)^{-2s+1}&\mbox{for}& r\ge a,
\end{array}
\right.
$$
where
\begin{equation}\label{eq:2.1}
\frac{1}{2}<s<\frac{\delta-2}{2}
\end{equation}
and $a=h^{-m}$ 
with some parameter $m>0$ to be fixed in the proof of Lemma 2.3 below depending only on $\delta$ and $s$.
Clearly, the first derivative (in sense of distributions) of $\mu$ satisfies
$$\mu'(r)=
\left\{
\begin{array}{lll}
 2(r+1)&\mbox{for}& 0\le r<a,\\
 (2s-1)(r+1)^{-2s}&\mbox{for}& r>a.
\end{array}
\right.
$$
The main properties of the functions $\mu$ and $\mu'$ are given in the following

\begin{lemma} For all $r>0$, $r\neq a$, we have the inequalities
\begin{equation}\label{eq:2.2}
2r^{-1}\mu(r)-\mu'(r)\ge 0,
\end{equation}
\begin{equation}\label{eq:2.3}
\mu'(r)\ge (2s-1)(r+1)^{-2s},
\end{equation}
\begin{equation}\label{eq:2.4}
\frac{\mu(r)^2}{\mu'(r)}\lesssim a^4(r+1)^{2s}.
\end{equation}
\end{lemma}

The proof of this lemma is straightforward and therefore we omit it. We now turn to the 
construction of the phase function
$\varphi\in C^1([0,+\infty))$ such that $\varphi(0)=0$ and $\varphi(r)>0$ for $r>0$. 
We define the first derivative of $\varphi$ by
$$\varphi'(r)=
\left\{
\begin{array}{lll}
 \tau(r+1)^{-1}- \tau(a+1)^{-1}&\mbox{for}& 0\le r\le a,\\
 0&\mbox{for}& r\ge a,
\end{array}
\right.
$$
where 
\begin{equation}\label{eq:2.5}
\tau=\tau_0h^{-1/3}
\end{equation}
with some parameter $\tau_0\gg 1$ independent of $h$ to be fixed in Lemma 2.3 below.
Clearly, the first derivative of $\varphi'$ satisfies
$$\varphi''(r)=
\left\{
\begin{array}{lll}
 -\tau(r+1)^{-2}&\mbox{for}& 0\le r<a,\\
 0&\mbox{for}& r>a.
\end{array}
\right.
$$

\begin{lemma} For all $r\ge 0$ we have the bound
\begin{equation}\label{eq:2.6}
h^{-1}\varphi(r)\lesssim h^{-4/3}\log\frac{1}{h}. 
\end{equation}
\end{lemma}

{\it Proof.} We have
$$\max\varphi=\tau\int_0^a\varphi'(r)dr\le \tau\int_0^a (r+1)^{-1}dr=\tau\log(a+1)$$
which clearly implies (\ref{eq:2.6}) in view of the choice of $\tau$ and $a$.
\eproof

For $r\neq a$, set
$$A(r)=\left(\mu\varphi'^2\right)'(r)$$
and
$$B(r)=\frac{\left(\mu(r)\left(h^{-1}(r+1)^{-\delta}+|\varphi''(r)|\right)\right)^2}{h^{-1}\varphi'(r)\mu(r)+\mu'(r)}.$$
The following lemma will play a crucial role in the proof of the Carleman estimates in the next section.

\begin{lemma} Given any $C>0$ independent of the variable $r$ and the parameters $h$, $\tau$ and $a$, there exist $\tau_0=\tau_0(C)>0$ and 
$h_0=h_0(C)>0$ so that for $\tau$ satisfying (\ref{eq:2.5}) and for all $0<h\le h_0$ we have the inequality
\begin{equation}\label{eq:2.7}
A(r)-CB(r)\ge -\frac{E}{2}\mu'(r)
\end{equation}
for all $r>0$, $r\neq a$.
\end{lemma}

{\it Proof.} For $r<a$ we have
$$A(r)=-\left(\varphi'^2\right)'(r)+\tau^2\partial_r\left(1-(r+1)(a+1)^{-1}\right)^2$$ 
$$=-2\varphi'(r)\varphi''(r)-2\tau^2(a+1)^{-1}\left(1-(r+1)(a+1)^{-1}\right)$$
$$\ge 2\tau(r+1)^{-2}\varphi'(r)-2\tau^2(a+1)^{-1}$$
$$\ge 2\tau(r+1)^{-2}\varphi'(r)-\tau^2 a^{-1}\mu'(r)$$
$$\ge 2\tau(r+1)^{-2}\varphi'(r)-{\cal O}(h^{m-1})\mu'(r)$$
where we have used that $\mu'(r)=2(r+1)$. Taking $m>2$ we get
\begin{equation}\label{eq:2.8}
A(r)\ge 2\tau(r+1)^{-2}\varphi'(r)-{\cal O}(h)\mu'(r)
\end{equation}
for all $r<a$. We will now bound the function $B$ from above. Let first $0<r\le \frac{a}{2}$. Since in this case we have
$$\varphi'(r)\ge \frac{\tau}{3}(r+1)^{-1}$$
 we obtain
$$B(r)\lesssim \frac{\mu(r)\left(h^{-2}(r+1)^{-2\delta}+\varphi''(r)^2\right)}{h^{-1}\varphi'(r)}$$ 
 $$\lesssim (\tau h)^{-1}\frac{\mu(r)(r+1)^{2-2\delta}}{\varphi'(r)^2}\tau(r+1)^{-2}\varphi'(r)
 + h\frac{\mu(r)\varphi''(r)^2}{\mu'(r)\varphi'(r)}\mu'(r)$$ 
$$\lesssim \tau^{-3}h^{-1}(r+1)^{6-2\delta}\tau(r+1)^{-2}\varphi'(r)+\tau h\mu'(r)$$
$$\lesssim \tau_0^{-3}\tau(r+1)^{-2}\varphi'(r)+ h^{2/3}\mu'(r)$$
where we have used that $\delta>3$. 
This bound together with (\ref{eq:2.8}) clearly imply (\ref{eq:2.7}), provided $\tau_0^{-1}$ and $h$ are taken small enough
depending on $C$.

Let now $\frac{a}{2}<r<a$. Then we have the bound
$$B(r)\le\left(\frac{\mu(r)}{\mu'(r)}\right)^2\left(h^{-1}(r+1)^{-\delta}+|\varphi''(r)|\right)^2\mu'(r)$$
$$\lesssim \left(h^{-2}(r+1)^{2-2\delta}+\tau^2(r+1)^{-2}\right)\mu'(r)$$
$$\lesssim \left(h^{-2}a^{2-2\delta}+\tau^2a^{-2}\right)\mu'(r)$$
$$\lesssim \left(h^{2m(\delta-1)-2}+h^{2m-2/3}\right)\mu'(r)\lesssim h\mu'(r)$$
provided $m$ is taken large enough. Again, this bound together with (\ref{eq:2.8}) imply (\ref{eq:2.7}).

It remains to consider the case $r>a$. Using that $\mu={\cal O}(a^2)$ together with (\ref{eq:2.3}) and taking into account
that $s$ satisfies (\ref{eq:2.1}), we get
$$B(r)=\frac{\left(\mu(r)\left(h^{-1}(r+1)^{-\delta}\right)\right)^2}{\mu'(r)}$$
$$\lesssim h^{-2}a^4(r+1)^{4s-2\delta}\mu'(r)\lesssim h^{-2}a^{4+4s-2\delta}\mu'(r)$$
$$\lesssim h^{2m(\delta-2-2s)-2}\mu'(r)\lesssim h\mu'(r)$$
provided that $m$ is taken large enough. Since in this case $A(r)=0$, the above bound clearly implies (\ref{eq:2.7}).
\eproof

\section{Carleman estimates} 

Our goal in this section is to prove the following

\begin{Theorem} Suppose (\ref{eq:1.1}) holds and let $s$ satisfy (\ref{eq:2.1}). Then, for all functions
$f\in H^2(\mathbb{R}^n)$ such that $\langle x\rangle^{s}(P(h)-E\pm i\varepsilon)f\in L^2$ and for all
$0<h\ll 1$, $0<\varepsilon\le ha^{-2}$, we have the estimate 
 \begin{equation}\label{eq:3.1}
\|\langle x\rangle^{-s}e^{\varphi/h}f\|_{L^2}\le Ca^2h^{-1}\|\langle x\rangle^{s}e^{\varphi/h}(P(h)-E\pm i\varepsilon)f\|_{L^2}
+Ca\tau(\varepsilon/h)^{1/2}\|e^{\varphi/h}f\|_{L^2}
\end{equation}
with a constant $C>0$ independent of $h$, $\varepsilon$ and $f$.
\end{Theorem}

{\it Proof.} We pass to the polar coordinates $(r,w)\in\mathbb{R}^+\times\mathbb{S}^{n-1}$, $r=|x|$, $w=x/|x|$, and
recall that $L^2(\mathbb{R}^n)=L^2(\mathbb{R}^+\times\mathbb{S}^{n-1}, r^{n-1}drdw)$. In what follows we denote by $\|\cdot\|$ and $\langle\cdot,\cdot\rangle$
the norm and the scalar product in $L^2(\mathbb{S}^{n-1})$. We will make use of the identity
\begin{equation}\label{eq:3.2}
 r^{(n-1)/2}\Delta  r^{-(n-1)/2}=\partial_r^2+\frac{\widetilde\Delta_w}{r^2}
\end{equation}
where $\widetilde\Delta_w=\Delta_w-\frac{1}{4}(n-1)(n-3)$ and $\Delta_w$ denotes the negative Laplace-Beltrami operator
on $\mathbb{S}^{n-1}$. Set $u=r^{(n-1)/2}e^{\varphi/h}f$ and
$${\cal P}^\pm(h)=r^{(n-1)/2}(P(h)-E\pm i\varepsilon) r^{-(n-1)/2},$$
$${\cal P}^\pm_\varphi(h)=e^{\varphi/h}{\cal P}^\pm(h)e^{-\varphi/h}.$$
Using (\ref{eq:3.2}) we can write the operator ${\cal P}^\pm(h)$ in the coordinates $(r,w)$ as follows
$${\cal P}^\pm(h)={\cal D}_r^2+\frac{\Lambda_w}{r^2}-E\pm i\varepsilon +V$$
where we have put ${\cal D}_r=-ih\partial_r$ and $\Lambda_w=-h^2\widetilde\Delta_w$. Since the function $\varphi$
depends only on the variable $r$, this implies
$${\cal P}^\pm_\varphi(h)={\cal D}_r^2+\frac{\Lambda_w}{r^2}-E\pm i\varepsilon -\varphi'^2+h\varphi''+
2i\varphi'{\cal D}_r+V.$$
For $r>0$, $r\neq a$, introduce the function
$$F(r)=-\langle (r^{-2}\Lambda_w-E-\varphi'(r)^2)u(r,\cdot),u(r,\cdot)\rangle+\|{\cal D}_ru(r,\cdot)\|^2$$
and observe that its first derivative is given by
$$F'(r)=\frac{2}{r}\langle r^{-2}\Lambda_wu(r,\cdot),u(r,\cdot)\rangle
+((\varphi')^2)'\|u(r,\cdot)\|^2$$
$$-2h^{-1}{\rm Im}\,\langle {\cal P}^\pm_\varphi(h)u(r,\cdot),{\cal D}_ru(r,\cdot)\rangle$$
$$\pm 2\varepsilon h^{-1}{\rm Re}\,\langle u(r,\cdot),{\cal D}_ru(r,\cdot)\rangle+4h^{-1}\varphi'\|{\cal D}_ru(r,\cdot)\|^2$$ 
$$+2h^{-1}{\rm Im}\,\langle (V+h\varphi'')u(r,\cdot),{\cal D}_ru(r,\cdot)\rangle.$$
Thus, if $\mu$ is the function defined in the previous section, we obtain the identity
$$\mu'F+\mu F'=
(2r^{-1}\mu-\mu')\langle r^{-2}\Lambda_wu(r,\cdot),u(r,\cdot)\rangle
+(E\mu'+(\mu(\varphi')^2)')\|u(r,\cdot)\|^2$$
$$-2h^{-1}\mu{\rm Im}\,\langle {\cal P}^\pm_\varphi(h)u(r,\cdot),{\cal D}_ru(r,\cdot)\rangle$$
$$\pm 2\varepsilon h^{-1}\mu{\rm Re}\,\langle u(r,\cdot),{\cal D}_ru(r,\cdot)\rangle+(\mu'+4h^{-1}\varphi'\mu)\|{\cal D}_ru(r,\cdot)\|^2$$ 
$$+2h^{-1}\mu{\rm Im}\,\langle (V+h\varphi'')u(r,\cdot),{\cal D}_ru(r,\cdot)\rangle.$$
Using that $\Lambda_w\ge 0$ together with (\ref{eq:2.2}) we get the inequality
$$\mu'F+\mu F'\ge (E\mu'+(\mu(\varphi')^2)')\|u(r,\cdot)\|^2+(\mu'+4h^{-1}\varphi'\mu)\|{\cal D}_ru(r,\cdot)\|^2$$
$$-\frac{3h^{-2}\mu^2}{\mu'}\|{\cal P}^\pm_\varphi(h)u(r,\cdot)\|^2-\frac{\mu'}{3}\|{\cal D}_ru(r,\cdot)\|^2$$
$$-\varepsilon h^{-1}\mu\left(\|u(r,\cdot)\|^2+\|{\cal D}_ru(r,\cdot)\|^2\right)$$
$$-3h^{-2}\mu^2(\mu'+4h^{-1}\varphi'\mu)^{-1}\|(V+h\varphi'')u(r,\cdot)\|^2
-\frac{1}{3}(\mu'+4h^{-1}\varphi'\mu)\|{\cal D}_ru(r,\cdot)\|^2$$
 $$\ge \left(E\mu'+(\mu(\varphi')^2)'-C\mu^2(\mu'+h^{-1}\varphi'\mu)^{-1}(h^{-1}(r+1)^{-\delta}
 +|\varphi''|)^2\right)\|u(r,\cdot)\|^2$$
$$-\frac{3h^{-2}\mu^2}{\mu'}\|{\cal P}^\pm_\varphi(h)u(r,\cdot)\|^2
-\varepsilon h^{-1}\mu\left(\|u(r,\cdot)\|^2+\|{\cal D}_ru(r,\cdot)\|^2\right)$$
with some constant $C>0$. Now we use Lemma 2.3 to conclude that
$$\mu'F+\mu F'\ge \frac{E}{2}\mu'\|u(r,\cdot)\|^2-\frac{3h^{-2}\mu^2}{\mu'}\|{\cal P}^\pm_\varphi(h)u(r,\cdot)\|^2$$ $$
-\varepsilon h^{-1}\mu\left(\|u(r,\cdot)\|^2+\|{\cal D}_ru(r,\cdot)\|^2\right).$$
We now integrate this inequality with respect to $r$ and use that, since $\mu(0)=0$, we have
$$\int_0^\infty(\mu'F+\mu F')dr=0.$$
Thus we obtain the estimate
\begin{equation}\label{eq:3.3}
\frac{E}{2}\int_0^\infty\mu'\|u(r,\cdot)\|^2dr\le 3h^{-2}\int_0^\infty\frac{\mu^2}{\mu'}\|{\cal P}^\pm_\varphi(h)u(r,\cdot)\|^2dr$$ $$
+\varepsilon h^{-1}\int_0^\infty\mu\left(\|u(r,\cdot)\|^2+\|{\cal D}_ru(r,\cdot)\|^2\right)dr.
\end{equation}
Using that $\mu={\cal O}(a^2)$ together with (\ref{eq:2.3}) and (\ref{eq:2.4}) we get from (\ref{eq:3.3})
\begin{equation}\label{eq:3.4}
\int_0^\infty(r+1)^{-2s}\|u(r,\cdot)\|^2dr\le Ca^4h^{-2}\int_0^\infty(r+1)^{2s}\|{\cal P}^\pm_\varphi(h)u(r,\cdot)\|^2dr$$ $$
+C\varepsilon h^{-1}a^2\int_0^\infty\left(\|u(r,\cdot)\|^2+\|{\cal D}_ru(r,\cdot)\|^2\right)dr
\end{equation}
with some constant $C>0$ independent of $h$ and $\varepsilon$. On the other hand, we have the identity
$${\rm Re}\,\int_0^\infty\langle {\cal P}^\pm_\varphi(h)u(r,\cdot),u(r,\cdot)\rangle dr=\int_0^\infty\|{\cal D}_ru(r,\cdot)\|^2dr
+\int_0^\infty \langle r^{-2}\Lambda_wu(r,\cdot),u(r,\cdot)\rangle dr$$
$$-\int_0^\infty(E+\varphi'^2)\|u(r,\cdot)\|^2dr+{\rm Re}\,\int_0^\infty\langle Vu(r,\cdot),u(r,\cdot)\rangle dr$$
which implies
\begin{equation}\label{eq:3.5}
\int_0^\infty\|{\cal D}_ru(r,\cdot)\|^2dr\le {\cal O}(\tau^2)\int_0^\infty\|u(r,\cdot)\|^2dr$$
$$+\gamma\int_0^\infty(r+1)^{-2s}\|u(r,\cdot)\|^2dr
+\gamma^{-1}\int_0^\infty(r+1)^{2s}\|{\cal P}^\pm_\varphi(h)u(r,\cdot)\|^2dr
\end{equation}
for every $\gamma>0$. We take now $\gamma$ small enough, independent of $h$, and recall that $\varepsilon h^{-1}a^2\le 1$. Thus, combining
the estimates (\ref{eq:3.4}) and (\ref{eq:3.5}), we get
\begin{equation}\label{eq:3.6}
\int_0^\infty(r+1)^{-2s}\|u(r,\cdot)\|^2dr\le Ca^4h^{-2}\int_0^\infty(r+1)^{2s}\|{\cal P}^\pm_\varphi(h)u(r,\cdot)\|^2dr$$ $$
+C\varepsilon h^{-1}a^2\tau^2\int_0^\infty\|u(r,\cdot)\|^2dr
\end{equation}
with a new constant $C>0$ independent of $h$ and $\varepsilon$. It is an easy observation now that the estimate 
(\ref{eq:3.6}) implies (\ref{eq:3.1}).
\eproof

\section{Resolvent estimates}

In this section we will derive the bound (\ref{eq:1.3}) from Theorem 3.1. Indeed, it follows from the estimate 
(\ref{eq:3.1}) and Lemma 2.2 that for $0<h\ll 1$, $0<\varepsilon\le ha^{-2}$ and $s$ satisfying (\ref{eq:2.1}) we have
\begin{equation}\label{eq:4.1}
\|\langle x\rangle^{-s}f\|_{L^2}\le M\|\langle x\rangle^{s}(P(h)-E\pm i\varepsilon)f\|_{L^2}
+M\varepsilon^{1/2}\|f\|_{L^2}
\end{equation}
where
$$M=\exp\left(Ch^{-4/3}\log(h^{-1})\right)$$
with a constant $C>0$ independent of $h$ and $\varepsilon$. On the other hand, since the operator $P(h)$ is symmetric, we have
\begin{equation}\label{eq:4.2}
\varepsilon\|f\|^2_{L^2}=\pm{\rm Im}\,\langle (P(h)-E\pm i\varepsilon)f,f\rangle_{L^2}$$
$$\le (2M)^{-2}\|\langle x\rangle^{-s}f\|^2_{L^2}+(2M)^2\|\langle x\rangle^{s}(P(h)-E\pm i\varepsilon)f\|^2_{L^2}.
\end{equation}
We rewrite (\ref{eq:4.2}) in the form
\begin{equation}\label{eq:4.3}
M\varepsilon^{1/2}\|f\|_{L^2}\le \frac{1}{2}\|\langle x\rangle^{-s}f\|_{L^2}+
2M^2\|\langle x\rangle^{s}(P(h)-E\pm i\varepsilon)f\|_{L^2}.
\end{equation}
We now combine (\ref{eq:4.1}) and (\ref{eq:4.3}) to get
\begin{equation}\label{eq:4.4}
\|\langle x\rangle^{-s}f\|_{L^2}\le 4M^2\|\langle x\rangle^{s}(P(h)-E\pm i\varepsilon)f\|_{L^2}.
\end{equation}
It follows from (\ref{eq:4.4}) that the resolvent estimate
\begin{equation}\label{eq:4.5}
\left\|\langle x\rangle^{-s}(P(h)-E\pm i\varepsilon)^{-1}\langle x\rangle^{-s}
\right\|_{L^2\to L^2}\le 4M^2
\end{equation}
holds for all $0<h\ll 1$, $0<\varepsilon\le ha^{-2}$ and $s$ satisfying (\ref{eq:2.1}).
On the other hand, for $\varepsilon\ge ha^{-2}$ the estimate (\ref{eq:4.5}) holds in a trivial way. Indeed, in this case, 
since the operator $P(h)$ is symmetric, the norm of the resolvent is upper bounded by $\varepsilon^{-1}=
{\cal O}(h^{-2m-1})$. Finally, observe that if (\ref{eq:4.5}) holds for $s$ satisfying (\ref{eq:2.1}),
it holds for all $s>1/2$.

\end{document}